\documentclass[reqno,10pt, centertags,draft]{amsart}
\usepackage{amsmath,amsthm,amscd,amssymb,latexsym,upref}

\makeatletter
\def\theequation{\@arabic\c@equation}

\newcommand{\bbN}{{\mathbb{N}}}
\newcommand{\bbR}{{\mathbb{R}}}

\newcommand{\bbC}{{\mathbb{C}}}

\newcommand{\cD}{{\mathcal D}}

\newcommand{\no}{\nonumber}
\newcommand{\lb}{\label}
\newcommand{\f}{\frac}

\newcommand{\ol}{\overline}

\newcommand{\Oh}{O}
\newcommand{\oh}{o}

\newcommand{\loc}{\text{\rm{loc}}}

\newcommand{\dom}{\text{\rm{dom}}}

\newcommand{\supp}{\text{\rm{supp}}}

\newcommand{\bi}{\bibitem}

\renewcommand{\ln}{\text{\rm ln}}

\numberwithin{equation}{section}

\newtheorem{theorem}{Theorem}[section]

\newtheorem{hypothesis}[theorem]{Hypothesis}
\theoremstyle{definition}
\newtheorem{definition}[theorem]{Definition}
\theoremstyle{remark}
\newtheorem{remark}[theorem]{Remark}
\newtheorem{example}[theorem]{Example}

\begin{document}

\title[Povzner--Wienholtz-type Results for Sturm--Liouville
Operators]{On Povzner--Wienholtz-type Self-Adjointness Results
for Matrix-Valued Sturm--Liouville Operators}
\author[S.\ Clark and F.\ Gesztesy]{Steve Clark and Fritz Gesztesy}
\address{Department of Mathematics and Statistics, University
of Missouri-Rolla, Rolla, MO 65409, USA}
\email{sclark@umr.edu}
\urladdr{http://www.umr.edu/\~{ }clark}
\address{Department of Mathematics,
University of Missouri, Columbia, MO 65211, USA}
\email{fritz@math.missouri.edu}
\urladdr{http://www.math.missouri.edu/people/fgesztesy.html}
\date{March 6, 2002.}
\thanks{Proc. Roy. Soc. Edinburgh {\bf A}, to appear.}

\begin{abstract}
We derive Povzner--Wienholtz-type self-adjointness results for $m\times m$
matrix-valued Sturm--Liouville operators
$T=R^{-1}\big[-\f{d}{dx}P\f{d}{dx}+Q\big]$ in $L^2((a,b);Rdx)^m$,
$m\in\bbN$, for $(a,b)$ a half-line or $\bbR$.
\end{abstract}

\maketitle

\section{Introduction} \lb{s1}

In two recent papers concerned with certain aspects of inverse
spectral theory for matrix-valued one dimensional Schr\"odinger
operators (cf. \cite{BGMS02}, \cite{GS02}), it became necessary to
infer the limit point case at $\pm\infty$ from nonoscillatory
properties of the Schr\"odinger differential expression at
$\pm\infty$. Put differently, it was necessary to conclude that
boundedness from below implies (essential) self-adjointness for
matrix-valued Schr\"odinger operators under minimal regularity assumptions
on the potential coefficient. While the problem of self-adjointness of
matrix-valued (and more generally, operator-valued) differential operators
has been studied in the literature (see, e.g.,
\cite{Li54}  and \cite{VG70}), apparently, no result assuming lower
semiboundedness of the associated minimal differential operator
alone, appears to be known. (Reference \cite{VG70}, e.g., assumes
additional local boundedness of the operator-valued potential in the
Schr\"odinger operator in question.) In this note we close this gap by
proving Povzner--Wienholtz-type self-adjointness results for
one-dimensional scalar and matrix-valued Sturm--Liouville operators under
optimal conditions on the potential coefficient.

Hartman~\cite{Ha48} and Rellich~\cite{Re51} were first to show
that boundedness below of the scalar Sturm--Liouville expression
$\ell = r^{-1}(-\frac{d}{dx}p\frac{d}{dx} +q)$ and  integral
conditions on the functions $r$ and $p$, which are necessarily
satisfied when $p,r=1$, are together sufficient to imply essential
self-adjointness of the minimal operator for $\ell$ (cf.
Section~\ref{s2}). Independently, Povzner~\cite{Po53} and
subsequently Wienholtz~\cite{Wi58} were first to show that
semi-boundedness for the operator $L_0=(-\Delta + q)\,
\upharpoonright\,C_0^{\infty}(\bbR^n)$ implies its essential
self-adjointness when $q$ is real-valued and continuous. This
result has since been shown for more general $q$ as may be seen by
the work and references found in \cite{Hi92}, \cite{Sc91},
\cite{Si78}, and \cite{Si92}. Closely related self-adjointness results
based on the concept of a finite rate of propagation were discussed by   
Berezanskii \cite[Sect.\ VI.1.7]{Be68} and further devloped, for
instance, in \cite{BS81}, \cite[Sect.\ 13.8]{BSU96}, \cite{Or88},
\cite{Ro85a}. These works typically address the scalar case for ordinary or
partial differential expressions, and not  the one-dimensional matrix
setting that is our interest. After completing our work we became aware of
a recent very detailed study of essential self-adjointness of general
Schr\"odinger operators of the type
$D^*D+V$, where $D$ is a first-order  elliptic differential operator
acting on the space of sections of a hermitian vector bundle over a
manifold, by Braverman, Milatovic, and Shubin \cite{BMS02}. This paper
appears to contain the state of the art in the multi-dimensional context. 

In Section~\ref{s2}, we define Sturm--Liouville operators for the
scalar case on the whole or half-line and recall in
Theorem~\ref{t2.3}  the equivalence between the boundedness below
of our operators and the nonoscillation of solutions of the
associated homogeneous equation near endpoints of the given
interval. The fundamental results of Hartman~\cite{Ha48} and of
Rellich~\cite{Re51} are recalled in Theorem~\ref{t2.4} for the
scalar Sturm--Liouville case. In Section~\ref{s3}, we obtain in
Theorem~\ref{t2.5} a result for the scalar Sturm--Liouville
operator, and in Theorem~\ref{t2.11} a result for the
matrix-valued Sturm--Liouville operator, like that of Hartman and
Rellich but independent of both as seen in Examples~\ref{e2.7} and
\ref{e2.8}.

Hartman proved his result using the concept of principal and
nonprincipal solutions while Rellich obtained a proof of this
result by focusing upon the boundedness below of the associated
minimal operator. Each of these approaches appears not to extend
to the matrix setting. Instead, our results are obtained using the
approach introduced by Wienholtz~\cite{Wi58} to prove the result
first due to Povzner~\cite{Po53}.
\section{Self-adjointness for the scalar Sturm--Liouville case}
\lb{s2}
Beginning with the scalar case, we introduce the following basic
assumptions.

\begin{hypothesis} \lb{h2.1} ${}$ \\
$(i)$ Let $c<b\leq\infty$ and suppose that $p,q,r$ are
$($Lebesgue$)$ measurable on $[c,b)$, and that
\begin{subequations}  \lb{2.1A}
\begin{align}
& p>0, r>0 \text{ a.e.\ on $[c,b)$, $q$ real-valued}, \lb{2.1aA} \\
& 1/p, q, r \in L^1 ([c,d]; dx) \text{ for all $d\in [c,b)$}. \lb{2.1bA}
\end{align}
\end{subequations}
$(ii)$ Let $-\infty\leq a<b\leq\infty$ and suppose that $p,q,r$ are
$($Lebesgue$)$ measurable on $(a,b)$, and that
\begin{subequations}  \lb{2.1}
\begin{align}
& p>0, r>0 \text{ a.e.\ on $(a,b)$, $q$ real-valued}, \lb{2.1a} \\
& 1/p, q, r \in L^1_{\loc} ((a,b); dx). \lb{2.1b}
\end{align}
\end{subequations}
\end{hypothesis}

The case $-\infty\leq a<c$ and $p,q,r$ (Lebesgue) measurable on $(a,c]$, 
$p>0, r>0$ a.e.\ on $(a,c]$, $q$ real-valued, and $1/p, q, r \in
L^1 ([d,c]; dx)$, for all $d\in (a,c]$, follows upon reflection from case
$(i)$ in Hypothesis \ref{h2.1} and hence will not be separately discussed
in the following.

Given Hypothesis~\ref{h2.1}, we consider the differential expression
\begin{equation}
\ell =\f{1}{r}\bigg(-\f{d}{dx}p\f{d}{dx}+q\bigg) \lb{2.2}
\end{equation}
on $[c,b)$ and $(a,b)$, respectively. Next, we introduce for
$c\in(a,b)$,
\begin{align}
\cD_{\max, \text{c,b}}&=\Big\{u\in L^2([c,b);rdx)\,\Big|\,
u, pu'\in AC([c,d]) \text{ for all $d\in [c,b)$;} \; \no \\
& \qquad u(c)=0; \,
 \ell u\in L^2([c,b);rdx)\Big\}, \lb{2.5} \\
\cD_{\min,\text{c,b}}&=\{u\in\cD_{\max,\text{c,b}} \,|\, \supp(u)\subset
[c,b) \text{ compact}\}, \lb{2.6} \\
\cD_{\max}&=\{u\in L^2((a,b);rdx)\,|\, u, pu'\in AC_{\loc} ((a,b)); 
\;\\
 &\qquad \ell u\in L^2((a,b);rdx)\}, \lb{2.3} \\
\cD_{\min}&=\{u\in\cD_{\max} \,|\, \supp(u)\subset (a,b) \text{
compact}\}, \lb{2.4} 
\end{align}
where $AC_{\loc}((a,b))$ denotes the set of locally absolutely continuous
functions on $(a,b)$ and $AC([c,d])$ represents the set of absolutely
continuous  functions on $[c,d]$.

Minimal operators $T_{\min,\text{c,b}}$ and $T_{\min}$ and maximal
operators $T_{\max,\text{c,b}}$ and $T_{\max}$ in $L^2((c,b);rdx)$ and 
$L^2((a,b);rdx)$, respectively, associated with $\ell$, are then defined by
\begin{align}
T_{\min,\text{c,b}}u&=\ell u, \quad u\in \dom
(T_{\min,\text{c,b}})=\cD_{\min,\text{c,b}}, \lb{2.11} \\
T_{\max,\text{c,b}}u&=\ell u, \quad u\in \dom
(T_{\max,\text{c,b}})=\cD_{\max,\text{c,b}}, \lb{2.12} \\
T_{\min}u&=\ell u, \quad u\in \dom (T_{\min})=\cD_{\min},  \lb{2.9} \\
T_{\max}u&=\ell u, \quad u\in \dom (T_{\max})=\cD_{\max},  \lb{2.10}
\end{align}
respectively. Then $T_{\min}$ and $T_{\min,\text{c,b}}$ are densely
defined and (cf.\ \cite[p.\ 64, 88]{Na68})
\begin{align}
{T_{\min,\text{c,b}}}^*&=T_{\max,\text{c,b}}, \lb{2.16} \\
{T_{\min}}^*&=T_{\max}. \lb{2.15}
\end{align}

In the following we will frequently refer to solutions $u(\cdot, z)$ of
$\ell u=z u$ for some $z\in\bbC$. Such  solutions are always assumed to
be distributional solutions, that is, we tacitly assume
\begin{equation}
u(\cdot,z), \, pu'(\cdot,z) \in AC([c,d]), \; d\in [c,b) \,
\text{ (resp., $AC_{\loc} ((a,b))$,
etc.)} \lb{2.18}
\end{equation}
in such a case.

Next we recall the standard notion of (non)oscillatory differential
expressions.

\begin{definition} \lb{d2.3}
Assume Hypothesis~\ref{h2.1}\,$(ii)$ and fix $c\in (a,b)$. Then $\ell$
is said to be {\it nonoscillatory near $a$ $($or $b$$)$} for some
$\lambda\in\bbR$ if and only if every real-valued solution $u$ of $\ell
u=\lambda u$ has finitely many zeros in $(a,c)$ (resp., $(c,b)$).
Otherwise, $\ell$ is called {\it oscillatory near $a$ $($resp., $b$$)$}.
\end{definition}

The following is a key result in this connection.

\begin{theorem} [\cite{Ha82}, \cite{Ka78}, \cite{Re51}, \cite{Ro85}, 
\cite{Sc01}] \lb{t2.3} ${}$ \\
Assume Hypothesis~\ref{h2.1}\,$(i)$. Then the following assertions are
equivalent: \\
$(i)$ $T_{\min,\text{c,b}}$ $($and hence any symmetric extension of
$T_{\min,\text{c,b}}$$)$ is bounded from below. \\
$(ii)$ There exists a $\lambda_0\in\bbR$ such that $\ell$ is
nonoscillatory near $b$ for all $\lambda < \lambda_0$. \\
Assume Hypothesis~\ref{h2.1}\,$(ii)$. Then the following assertions are
equivalent: \\
$(iii)$ $T_{\min}$ $($and hence any symmetric extension of $T_{\min}$$)$
is bounded from below. \\
$(iv)$ There exists a $\lambda_0\in\bbR$ such that $\ell$ is
nonoscillatory near $a$ and $b$ for all $\lambda < \lambda_0$.
\end{theorem}

Next, we mention a fundamental result which links the nonoscillatory
behavior at one end point with the limit point property at that endpoint.

\begin{theorem} [Hartman \cite{Ha48} (cf.\ also \cite{Ge93}, \cite{Re51})]
\lb{t2.4}
Suppose Hypothesis~\ref{h2.1}\,$(i)$ and assume that for some
$\lambda_0\in\bbR$, $\ell-\lambda_0$ is nonoscillatory near $b$, or
equivalently, that $T_{\min,\text{c,b}}$ is bounded from below. Then, if
\begin{equation}
\Bigg|\int^b dx\,\bigg(\f{r(x)}{p(x)}\bigg)^{1/2}\Bigg| =\infty, \lb{2.19}
\end{equation}
$\ell$ is in the limit point case at $b$.
\end{theorem}

Originally, Hartman proved the special case $r=1$ by an elegant
application of the concept of (non)principal solutions $($the proof easily
extends to the general situation $r\neq 1$ described in 
Theorem~\ref{t2.4}$)$. Theorem~\ref{t2.4} was subsequently derived by
Rellich \cite{Re51} by focusing on operators  $T_{\min,\text{c,b}}$ 
bounded from below. (Actually, Hartman and Rellich assume $p,p'\in
C((a,b))$ and $q,r$ piecewise continuous in $(a,b)$ in addition to 
Hypothesis~\ref{h2.1}\,(i) but this is easily seen to be unnecessary.)

For an extension of Theorem~\ref{t2.4} to polynomials of $\ell$ we refer
to Read \cite{Re76}.

For essential self-adjointness results using conditions of the type 
\eqref{2.19} but with semiboundedness from below of the corresponding
minimal differential operators replaced by alternative conditions, we
refer to the recent work by Lesch and Malamud \cite{LM02}. The latter
treats the case of matrix-valued Schr\"odinger operators as well as the
case of general canonical systems. 

Finally, for extensions of Theorem~\ref{t2.4} to the multi-dimensional
case we refer to our account in the introduction. 
\section{Povzner--Wienholtz-type self-adjointness results} \lb{s3}

Next we formulate a result that resembles one by Povzner \cite{Po53} in
the context of partial differential operators. We will adapt a method of
proof due to Wienholtz \cite{Wi58}, who independently proved
Povzner's result. (Wienholtz's proof is also reproduced in
Glazman \cite[Chapter1, Theorem\ 35]{Gl65}.)

\begin{theorem} \lb{t2.5}
Assume Hypothesis~\ref{h2.1}\,$(i)$ with $b=\infty$ and $p\in
AC([c,c+\rho])$ for all $\rho>0$ and $r^{-1}(p')^2, r^{-1}\in
L^1_{\loc}((c,\infty))$. Moreover, suppose that
$T_{\min,\text{c},\infty}$ is bounded from below. Then, if
\begin{equation}
\|p/r\|_{L^\infty((\rho/2,\rho))}
\underset{\rho\uparrow\infty}{=}\Oh(\rho^2), \lb{2.22}
\end{equation}
$\ell$ is in the limit point case at $\infty$.
\end{theorem}
\begin{proof}
By von Neumann's theory of self-adjoint extensions (cf.\ \cite[p.\
137]{RS75}), to prove that $\ell$ is in the limit point case at
$\infty$, or equivalently, that $T_{\min,\text{c},\infty}$ is essentially
self-adjoint and hence its closure $T_{\max,\text{c},\infty}$ is
self-adjoint (since
$\ol{T_{\min,\text{c},\infty}}=T_{\min,\text{c},\infty}^*
=T_{\max,\text{c},\infty}$), we need to show that
\begin{equation}
\ker(T_{\min,\text{c},\infty}^*-C)=\{0\} \lb{2.23}
\end{equation}
for $C>0$ sufficiently large. Next, we assume, without loss of
generality, that for $u\in \dom (T_{\min,\text{c},\infty})$
\begin{equation}\lb{2.24}
\int_c^\infty r(x)dx\, \ol{u(x)}(\ell u)(x) > \int_c^\infty r(x)dx
\,|u(x)|^2 ,
\end{equation}
so we may take $C=0$ in \eqref{2.23}. In addition, we choose $c=0$
for simplicity.

To start the proof we assume that $\ell$ is not in the limit point case
at $\infty$. Then, by Weyl's alternative, all solutions $u$ of
\begin{equation}\lb{2.25}
\ell u =0
\end{equation}
are necessarily in $L^2([0,\infty);rdx)$, and hence there exists a
nontrivial real-valued solution $\hat u\in L^2([0,\infty);rdx)$ of
\eqref{2.25} that satisfies the Dirichlet boundary condition
$\hat u(0)=0$. To complete the proof, it is sufficient to show that
\eqref{2.24} rules out the existence of such a solution 
$\hat u\in L^2([0,\infty);rdx)$ of \eqref{2.25}.

Let
\begin{equation}
u_\rho(x) = \theta_\rho(x) \hat u(x)\, , \quad 0\le x <\infty, \lb{2.26}
\end{equation}
where
\begin{equation}
\theta_\rho(x) =\theta(x/\rho), \quad \theta \in
\bbC^\infty([0,\infty)), \quad \theta(x) =\begin{cases}
 1 & 0 \le x \le 1/2, \\ 0 & x >1. \end{cases} \lb{2.27}
\end{equation}
As defined, $u_\rho\in\dom (T_{\min,0,\infty})$. Indeed, since
$p\in AC([0,\rho])$ for all $\rho>0$ by hypothesis, one verifies that
\begin{align}
\ell u_\rho &=\ell (\theta_\rho \hat u) =\theta_\rho(\ell \hat u)
-r^{-1}(p\theta_\rho')'\hat u-2r^{-1}\theta_\rho'(p\hat u') \no \\
&=-r^{-1}p'\theta_\rho'\hat u - r^{-1}p\theta_\rho''\hat u
-2r^{-1}\theta_\rho'(p\hat u'). \lb{2.28}
\end{align}
Since 
\begin{align}
&\int_{\rho/2}^\rho dx\,r^{-1}(p')^2(\theta_\rho')^2\hat u^2 \leq 
\|(\theta_\rho')^2\hat u^2\|_{L^\infty((\rho/2,\rho))} 
\int_{\rho/2}^\rho dx\,r^{-1}(p')^2 <\infty, \lb{2.28a} \\
&\int_{\rho/2}^\rho dx\, (r^{-1} p)\big[p(\theta_\rho'')^2\hat u^2\big]
\leq  (\rho/2)\|p/r\|_{L^\infty((\rho/2,\rho))} 
\|p\|_{L^\infty((\rho/2,\rho))} \no \\
& \hspace*{4.2cm} \times \|(\theta_\rho'')^2 
\hat u^2\|_{L^\infty((\rho/2,\rho))} <\infty, \lb{2.28b} \\
\intertext{and}
& \int_{\rho/2}^\rho dx\,r^{-1} (\theta_\rho')^2 (p\hat u')^2 \leq 
\|(\theta_\rho')^2\|_{L^\infty((\rho/2,\rho))} 
\|(p\hat u')^2\|_{L^\infty((\rho/2,\rho))} \no \\
& \hspace*{3.8cm} \times \int_{\rho/2}^\rho dx\,r^{-1} <\infty, \lb{2.28c}
\end{align}
one obtains $u_\rho\in\dom (T_{\min,0,\infty})$. Hence $u_\rho$ satisfies
\eqref{2.24}, that is,
\begin{equation}\lb{2.29}
\int_0^\infty r(x)dx\, \theta_\rho(x)^2 \hat u(x)^2 < \int_0^\infty
r(x)dx\, u_\rho(x)(\ell u_\rho)(x).
\end{equation}
Since 
\begin{equation}
ru_\rho (\ell u_\rho)=-\theta_\rho(p\theta_\rho')'\hat u^2 - 2\theta_\rho
\theta_\rho'\hat u (p\hat u')=(\theta_\rho')^2p\hat u^2 
-\big[\theta_\rho\theta_\rho' p\hat u^2 \big]' \lb{2.30}
\end{equation}
and $\theta_\rho(\rho)=\theta_\rho'(0)=0$, one obtains
\begin{align}
& \int_0^{\rho/2} r(x)dx\, \hat u(x)^2 
 <\int_0^\infty r(x)dx\,\theta_\rho(x)^2 \hat u(x)^2 \no \\
& < \int_0^\infty r(x)dx\, u_\rho(x)(\ell u_\rho)(x)  
=\int_0^\infty dx\, \theta_\rho'(x)^2 p(x)\hat u(x)^2 \no \\
& =\int_0^\infty r(x)dx\, \theta_\rho'(x)^2 [p(x)/r(x)] 
\hat u(x)^2 \no \\
& \leq \frac{C^2}{\rho^2} \|p/r\|_{L^\infty((\rho/2,\rho))} \int_{\rho/2}^\rho
r(x)dx\, \hat u(x)^2  \no \\
&\underset{\rho\uparrow\infty}{=} \Oh(1)  \int_{\rho/2}^\rho r(x)dx\, \hat
u(x)^2 \underset{\rho \uparrow\infty}{=} \oh(1) \lb{2.31}
\end{align}
for some constant $C>0$, implying $\hat u=0$ a.e. on $[0,\infty)$.
\end{proof}

\begin{remark} \lb{r2.6}
Even though Theorems~\ref{t2.4} (for $b=\infty$) and \ref{t2.5} show some
resemblance,  their proofs follow quite different strategies. In
particular, the proof of Theorem~\ref{t2.5} seems to work only for
half-lines (i.e., $b=\infty$). Moreover, as the following examples show,
neither theorem implies the other.
\end{remark}

\begin{example} \lb{e2.7}
Consider $c=0$, $b=\infty$, $q=0$, $p,r\in C^\infty([0,\infty))$, and
\begin{equation}
\f{r(x)}{p(x)}=\f{1}{x^2\ln(x)} \, \text{ for $x\geq e$}. \lb{2.36}
\end{equation}
Then one easily verifies that condition \eqref{2.19} is satisfied while
\eqref{2.22} is not.
\end{example}

\begin{example} \lb{e2.8}
Consider $c=0$, $b=\infty$ $q=0$, $p,r\in C([0,\infty))$, and
\begin{equation}
\f{r(x)}{p(x)}=\begin{cases} \f{1}{x^2}, & x\in [n,n+(1/n)], \; n\in\bbN,
\\ 0, & \text{otherwise}. \end{cases} \lb{2.37}
\end{equation}
Then one easily verifies that condition \eqref{2.22} is satisfied while
\eqref{2.19} is not.
\end{example}

We note that $q=0$ implies the nonoscillatory behavior of $\ell$ since
solving $\ell u=0$ readily yields
\begin{equation}
u(x)=C\int^x dx'\,p(x')^{-1} + D, \quad C,D \in\bbC, \lb{2.38}
\end{equation}
and clearly $u$ in \eqref{2.38} has at most one zero. Moreover, if $r$ is
not integrable at $\infty$, there exists a particular solution
$u_0=D\notin L^2([0,\infty);rdx)$ (choosing $D\neq 0$ in \eqref{2.38}) of
$\ell u=0$ and hence $\ell$ is also in the limit point case in this
situation.

\begin{remark} \lb{r2.9}
Examining the proof of Theorem~\ref{t2.5}, it is possible to relax the
assumption $p\in AC_{\loc}(\bbR)$ by an appropriate modification of
condition \eqref{2.22}. In fact, modifying the cutoff function
$\theta_\rho$ in \eqref{2.27} as follows,
\begin{align}
&\theta_\rho\in C^\infty ([0,\infty)), \quad
\theta_\rho(x)=\begin{cases} 1, & 0\leq x\leq \rho/2, \\
f_\rho\big(\rho^{-1}\int_0^x dx' \, p(x')^{-1}\big), & \rho/2\leq x \leq
\rho, \\ 0, &x\geq \rho, \end{cases} \lb{2.39} \\
&f_\rho\Big(\rho^{-1}\smallint_0^{\rho/2}dx\, p(x)^{-1}\Big)=1, \quad
f_\rho\Big(\rho^{-1}\smallint_0^{\rho}dx\,p(x)^{-1}\Big)=0 \lb{2.40}
\end{align}
for an appropriate smooth function $f_\rho$, and assuming
\begin{equation}
\bigg\|(pr)^{-1}\Big[f'_\rho\Big(\rho^{-1}\smallint_0^{\cdot}dx'\,
p(x')^{-1}\Big)\Big]^2\bigg\|_{L^\infty((\rho/2,\rho))}
\underset{\rho\uparrow\infty}{=}\Oh(\rho^2), \lb{2.41}
\end{equation}
instead of condition \eqref{2.22}, one verifies that
\begin{align}
\ell u_\rho &=\ell (\theta_\rho \hat u) =\theta_\rho(\ell \hat u)
-r^{-1}(p\theta_\rho')'\hat u-2r^{-1}\theta_\rho'(p\hat u') \no \\
&=-r^{-1}(p\theta_\rho')'\hat u -2r^{-1}\theta_\rho'(p\hat u') \no \\
&=-\rho^{-2}(pr)^{-1}
f_\rho''\Big(\rho^{-1}\smallint_0^xdx'\,p(x')^{-1}\Big)\hat u \no \\
& \quad -2\rho^{-1}(pr)^{-1}f_\rho'\Big(\rho^{-1}
\smallint_0^xdx'\,p(x')^{-1}\Big)(p\hat u')  \lb{2.42}
\end{align}
is well-defined. Following the proof of Theorem~\ref{t2.5} step by step
then results in
\begin{align}
&\int_0^\infty dx\, \theta_\rho'(x)^2p(x)\hat u(x)^2
=\int_0^\infty r(x)dx\, \theta_\rho'(x)^2 p(x)r(x)^{-1} \hat u(x)^2 \no \\
&\leq\frac{C^2}{\rho^2}\bigg\|(pr)^{-1}\Big[f'_\rho\Big(\rho^{-1}
\smallint_0^{\cdot}dx'\,
p(x')^{-1}\Big)\Big]^2\bigg\|_{L^\infty((\rho/2,\rho))}
\int_{\rho/2}^\rho r(x)dx\, \hat u(x)^2 \lb{2.43}
\end{align}
for some constant $C>0$ and hence again in $\hat u=0$ a.e.\ on
$[0,\infty)$ as in \eqref{2.31}.
\end{remark}

When trying to extend Theorem~\ref{t2.4} to the matrix-valued context,
one observes that the use of (matrix-valued) principal solutions can
indeed exclude the limit circle case at $b$ under the assumption of
$T_{\min,\text{c,b}}$ being bounded from below. However, it is not clear
how to exclude all intermediate cases between those of limit circle and
limit point that we may then infer the limit point case at $b$. The
strategy of proof of Theorem~\ref{t2.5}, on the other hand, will now be
shown to extend to the matrix-valued case in due course.

Throughout the rest of this note all matrices will be considered over the
field of complex numbers $\bbC$, and the corresponding linear space of
$k\times\ell$ matrices will be denoted by $\bbC^{k\times\ell}$,
$k,\ell\in\bbN$.  Moreover, $I_k$ denotes the identity matrix in
$\bbC^{k\times k}$ and $M^*$ the adjoint (i.e., complex conjugate
transpose) of a matrix $M$. A positive definite matrix $P$ is
denoted by $P>0$, a nonnegative matrix $Q$ by $Q\geq 0$.

We start by introducing the following hypothesis.

\begin{hypothesis} \lb{h2.10} ${}$ \\
Let $c<b\leq\infty$ and suppose that $P,Q,R\in\bbC^{m\times m}$ are
$($Lebesgue$)$ measurable on $[c,b)$, and that
\begin{subequations}  \lb{2.45A}
\begin{align}
& P>0, R>0 \text{ a.e.\ on $[c,b)$, $Q=Q^*$ self-adjoint}, \lb{2.45a} \\
& P^{-1}, Q, R \in L^1 ([c,d]; dx)^{m\times m} \text{ for all
$d\in [c,b)$}. \lb{2.45b}
\end{align}
\end{subequations}
\end{hypothesis}

As in the scalar context, the analogous case of the half-line $(a,c]$
follows upon reflection from the case $[c,b)$ in Hypothesis \ref{h2.10}
and hence is not separately discussed in the following.

Given Hypothesis~\ref{h2.10}, we consider the $m\times m$
matrix-valued differential expression
\begin{equation}
L = R^{-1}\bigg(-\f{d}{dx}P\f{d}{dx}+Q\bigg) \lb{2.47}
\end{equation}
on $[c,b)$. Next, we introduce 
\begin{align}
\cD_{\max, \text{c,b}}&=\Big\{u\in L^2([c,b);Rdx)^m\,\Big|\,
u, Pu'\in AC([c,d])^m \text{ for all $d\in [c,b)$;} \; \no \\
& \qquad u(c)=0; \,
 \ell u\in L^2([c,b);Rdx)^m\Big\}, \lb{2.48} \\
\cD_{\min,\text{c,b}}&=\{u\in\cD_{\max,\text{c,b}} \,|\, \supp(u)\subset
[c,b) \text{ compact}\}, \lb{2.49} 
\end{align}
where $AC([c,d])^m$ denotes the set of absolutely continuous vectors in
$\bbC^m$ on $[c,d]$. Moreover, $L^2(\Omega;Rdx)^m$ denotes the Hilbert
space of $\bbC^m$-valued (Lebesgue) measurable vectors $u,v$ on
$\Omega\subseteq\bbR$ corresponding to the $m\times m$
matrix-valued measure $Rdx$ with scalar product
$(u,v)_{L^2(\Omega;Rdx)^m}$, given by
\begin{equation}
(u,v)_{L^2(\Omega;Rdx)^m}=\int_\Omega dx \,
(u(x),R(x)v(x))_{\bbC^m}, \lb{2.52}
\end{equation}
where $ (f,g)_{\bbC^m}=\sum_{k=1}^m \ol{f}_kg_k, \; f,g
\in\bbC^m$. For notational simplicity we will abbreviate the
scalar product in ${L^2(\Omega;Rdx)^m}$ by $(u,v)$ only.

Minimal operators $T_{\min,\text{c,b}}$ and maximal operators
$T_{\max,\text{c,b}}$ in $L^2((c,b);Rdx)^m$ associated with $L$ are then
defined as in the scalar context by
\begin{align}
T_{\min,\text{c,b}}u&=L u, \quad u\in \dom
(T_{\min,\text{c,b}})=\cD_{\min,\text{c,b}}, \lb{2.53} \\
T_{\max,\text{c,b}}u&=L u, \quad u\in \dom
(T_{\max,\text{c,b}})=\cD_{\max,\text{c,b}}, \lb{2.54} 
\end{align}
respectively. Then $T_{\min}$ and $T_{\min,\text{c,b}}$ are densely
defined and \eqref{2.16}--\eqref{2.15} continue to hold in the
matrix-valued context.

After these preparations we are ready for the matrix-valued
extension of Theorem~\ref{t2.5}, the principal objective of this note. (We
denote by $\|M\|$ the operator norm of an $m\times m$ matrix $M$.)

\begin{theorem} \lb{t2.11}
Assume Hypothesis~\ref{h2.10} with $b=\infty$ and $P\in
AC([c,c+\rho])^{m\times m}$ for all $\rho>0$ and $\big\|P'R^{-1}P'\big\|,
\big\|R^{-1}\big\| \in L^1_{\loc} ((c,\infty))$. Moreover, suppose that
$T_{\min,\text{c},\infty}$ is bounded from below. Then, if
\begin{equation}
\text{\rm{ess\,sup}}_{x\in
[\rho/2,\rho]} \big\|R(x)^{-1/2}P(x)R(x)^{-1/2}\big\|
\underset{\rho\uparrow\infty}{=}\Oh(\rho^2), \lb{2.57} 
\end{equation}
$L$ is in the limit point case at $\infty$.
\end{theorem}
\begin{proof}
The strategy of proof is completely analogous to the scalar context,
Theorem~\ref{t2.5}. As before we need to show that
\begin{equation}
\ker(T_{\min,\text{c},\infty}^*-C)=\{0\} \lb{2.58}
\end{equation}
for $C>0$ sufficiently large. Again we assume, without loss of
generality, that for $u\in \dom (T_{\min,\text{c},\infty})$
\begin{equation}\lb{2.59}
\int_c^\infty dx\, (u(x),R(x)(L u)(x))_{\bbC^m} > \int_c^\infty
dx \,(u(x),R(x)u(x))_{\bbC^m},
\end{equation}
so we may take $C=0$ in \eqref{2.58}. Moreover, we choose $c=0$ as before.

To start the proof we assume that $L$ is not in the limit point case
at $\infty$. Then, by the matrix analog of Weyl's alternative (cf., e.g., 
\cite{KR74}, \cite{LM02}), $m+1$ solutions $u$ of
\begin{equation}\lb{2.60}
L u =0
\end{equation}
are necessarily in $L^2([0,\infty);Rdx)^m$, and hence there exists a
nontrivial real-valued solution $\hat u\in L^2([0,\infty);Rdx)$ of
\eqref{2.25} that satisfies the Dirichlet boundary condition $\hat u(0)=0$.
To complete the proof, it is sufficient to show that \eqref{2.59} rules
out the existence of such a solution $\hat u$.

Let
\begin{equation}
u_\rho(x) = \theta_\rho(x) \hat u(x)\, , \quad 0\le x <\infty, \lb{2.61}
\end{equation}
where
\begin{equation}
\theta_\rho(x) =\theta(x/\rho), \quad \theta \in
\bbC^\infty([0,\infty)), \quad \theta(x) =\begin{cases}
 1 & 0 \le x \le 1/2, \\ 0 & x >1. \end{cases} \lb{2.62}
\end{equation}
As defined, $u_\rho\in\dom (T_{\min,0,\infty})$. Indeed, since
$P\in AC([0,\rho])$ for all $\rho>0$ by hypothesis, one verifies that
\begin{align}
L u_\rho &=L (\theta_\rho \hat u) =\theta_\rho(L \hat u)
-R^{-1}(P\theta_\rho')'\hat u-2R^{-1}\theta_\rho'(P\hat u') \no \\
&=-R^{-1}P'\theta_\rho'\hat u -R^{-1}P\theta_\rho''\hat u
-2R^{-1}\theta_\rho'(P\hat u'). \lb{2.63}
\end{align}
Since 
\begin{align}
& \int_{\rho/2}^\rho dx\,\big(\theta_\rho'\hat u,P'R^{-1}P'\theta_\rho'\hat
u\big)_{\bbC^m} \leq \text{ess\,sup}_{x\in [\rho/2,\rho]}\|\theta_\rho'(x) 
\hat u(x)\|^2_{\bbC^m} \no \\
& \hspace*{4.9cm} \times \int_{\rho/2}^\rho dx\, \big\|P'R^{-1}P'\big\|
<\infty, \lb{2.63a} \\
& \int_{\rho/2}^\rho dx\,\big(P\theta_\rho''\hat u,R^{-1}P\theta_\rho''\hat
u\big)_{\bbC^m}
=\int_{\rho/2}^\rho dx\,\big(P^{1/2}\theta_\rho''\hat u,
P^{1/2}R^{-1}P^{1/2}P^{1/2}\theta_\rho''\hat u\big)_{\bbC^m} \no \\
& \leq (\rho/2)\text{ess\,sup}_{x\in [\rho/2,\rho]}\big\|P(x)^{1/2}
\big\|^2 \, \text{ess\,sup}_{x\in [\rho/2,\rho]}\|\theta_\rho''(x) \hat
u(x)\|^2_{\bbC^m} \no \\
& \quad \times \text{ess\,sup}_{x\in [\rho/2,\rho]}
\big\|P(x)^{1/2}R(x)^{-1}P(x)^{1/2}\big\| \no \\ 
&= (\rho/2)\text{ess\,sup}_{x\in [\rho/2,\rho]}\big\|P(x)^{1/2}\big\|^2 \,
\text{ess\,sup}_{x\in [\rho/2,\rho]}\|\theta_\rho''(x) \hat
u(x)\|^2_{\bbC^m} \no \\
& \quad \times \text{ess\,sup}_{x\in [(\rho/2,\rho)]} 
\big\|R(x)^{-1/2}P(x)R(x)^{-1/2}\big\| < \infty, \lb{2.63b} \\
\intertext{and}
& \int_{\rho/2}^\rho \big(R^{-1}\theta_\rho'(P\hat u'),RR^{-1}\theta_\rho'
(P\hat u')\big)_{\bbC^m} \no \\
& \leq \big\|(\theta_\rho')^2\big\|_{L^\infty([\rho/2,\rho])} \,
\text{ess\,sup}_{x\in [\rho/2,\rho]}\big\|P(x) \hat u'(x)\big\|^2_{\bbC^m} 
\int_{\rho/2}^\rho dx\, \big\|R^{-1}\big\| < \infty, \lb{2.63c}
\end{align}
and thus $u_\rho\in\dom (T_{\min,0,\infty})$. Hence $u_\rho$ satisfies
\eqref{2.59}, that is,
\begin{equation}\lb{2.64}
\int_0^\infty dx\, \theta_\rho(x)^2 (\hat u(x),R(x)\hat u(x))_{\bbC^m} <
\int_0^\infty dx\, (u_\rho(x),R(x)(L u_\rho)(x))_{\bbC^m}.
\end{equation}
Since
\begin{align}
\big[(u\rho,R(Lu_\rho))_{\bbC^m}+(R(Lu_\rho),u_\rho)_{\bbC^m} \big]/2 
=\big(\theta_\rho'\big)^2 (\hat u,P\hat u)_{\bbC^m} 
-\big[\theta_\rho\theta_\rho' (\hat u,P\hat u)_{\bbC^m}\big]' \lb{2.65}
\end{align}
and $\theta_\rho(\rho)=\theta_\rho'(0)=0$, one obtains
\begin{align}
& \int_0^{\rho/2} dx\, (\hat u(x),R(x)\hat u(x))_{\bbC^m}  <\int_0^\infty
dx\,\theta_\rho(x)^2 (\hat u(x),R(x)\hat u(x))_{\bbC^m} \no \\
& < \int_0^\infty (u_\rho(x),R(x)(L u_\rho)(x))_{\bbC^m} \no \\
& =(1/2) \int_0^\infty (u_\rho(x),R(x)(L u_\rho)(x))_{\bbC^m} + 
(1/2) \int_0^\infty (R(x)(L u_\rho)(x),u_\rho(x))_{\bbC^m} \no \\
& =\int_0^\infty \theta_\rho'(x)^2 (\hat u(x),P(x) \hat u(x))_{\bbC^m} 
\no \\ 
& =\int_0^\infty dx\, \theta_\rho'(x)^2
\big(R(x)^{1/2}\hat u(x),
R(x)^{-1/2}P(x)R(x)^{-1/2}R(x)^{1/2}\hat u(x)\big)_{\bbC^m} \no \\
&\leq\frac{C^2}{\rho^2} \text{ess\,sup}_{x\in [\rho/2,\rho]}
\big\|R(x)^{-1/2}P(x)R(x)^{-1/2}\big\|
\int_{\rho/2}^\rho dx\, (\hat u(x),R(x)\hat u(x))_{\bbC^m} \no \\
& \underset{\rho\uparrow\infty}{=} \Oh(1)  \int_{\rho/2}^\rho dx\,
(\hat u(x),R(x)\hat u(x))_{\bbC^m}
\underset{\rho \uparrow\infty}{=} \oh(1) \lb{2.66}
\end{align}
for some constant $C>0$, implying $\hat u=0$ a.e. on $[0,\infty)$.
\end{proof}

\bigskip\bigskip
\noindent {\bf Acknowledgements.} We are indebted to Mark Malamud and
Fedor  Rofe-Beketov for useful comments and hints to the literature. 
We also thank the referee for a very careful reading of our manuscript and
for pointing out a number of improvements concerning the presentation of
our results. 


\end{document}